\documentclass{article}
\usepackage{row}
\usepackage{graphicx}
\begin{document}
\title{Geometry in the Age of Enlightenment}
\author{Raymond O. Wells, Jr.
\footnote{Jacobs University Bremen; University of Colorado at Boulder; row@raw.com}}

\maketitle
\tableofcontents

\section{Introduction}
The {\it Age of Enlightenment} is a term that refers to a time of dramatic changes in western society in the arts, in science, in political thinking, and, in particular, in philosophical discourse.  It is generally recognized as being the period from the mid 17th century to the latter part of the 18th century.  It was a successor to the renaissance and reformation periods and was followed by what is termed the romanticism of the 19th century. In his book {\it A History of Western Philosophy} \cite{russell1946} Bertrand Russell (1872--1970) gives a very lucid description of this time period in intellectual history, especially in Book III, Chapter VI--Chapter XVII. He singles out Ren\'{e} Descartes as being the founder of the era of new philosophy in 1637 and continues to describe other philosophers who also made significant contributions to mathematics as well, such as Newton and Leibniz. This time of intellectual fervor included literature (e.g. Voltaire), music and the world of visual arts as well. One of the most significant developments was perhaps in the political world: here the absolutism of the church and of the monarchies were  questioned by the political philosophers of this era, ushering in the Glorious Revolution in England (1689), the American Revolution (1776), and the bloody French Revolution (1789).  All of these were culminations of this Age of Enlightenment which  permanently changed the shape of Western Civilization from the absolutism of the Middle Ages. A most important development of this time was the rise of science as it began to play an increasingly important role in the world at large, along with the technological advances which accompanied it.  Russell describes this advance in science very succinctly in his book.  

Mathematics experienced, as a part of this intellectual development, exciting growth with numerous new sets of ideas.  In this paper we would like to outline some of the important developments and new  ideas in geometry which were a part of this era.  There were, of course many other important mathematical development of this period such as in analysis, number theory, algebra, applied mathematics and many others.  We have described in two other papers \cite{wells2013} and \cite{wells2015} the very innovative geometric ideas that arose in the  fertile 19th century which extended many aspects of geometry to abstractions of various kinds. These ideas extended beyond the usual  study of geometry in two- and three-dimensional space as was inspired  by the mathematicians of the Greek era. However, in the Age of Enlightenment there were crucial new discoveries concerning curves and surfaces in the plane and in ordinary three space which became crucial building blocks for 19th century geometry. 

Between the time of the Greeks (which lasted about one millennium from {\it c.} 600 BC till {\it c.} 400 AD) and the rise of mathematical thought in Western Europe, in the 17th and 18th centuries, there were numerous mathematical developments in the Arabic and Indian cultures, primarily in arithmetic and algebra.  One of the most significant accomplishments of the Arabic world after the fall of the Roman empire and preceding the time of the Renaissance was the actual preservation of a substantive amount of the accomplishments of the Greek mathematicians.  In the realm of geometry, and in particular in the theory of curves and surfaces in a three-dimensional Euclidean space, there are almost no discernible developments during this interim time period. Pappus, whose work plays an important role in this context, came towards the very end of Greek mathematical culture in the fourth century A.D.. For more than 1000 years, mainly from the time of Pappus to the extraordinarily original work of Descartes in the 17th century, geometry seemed to be at a standstill.  

There are four major areas of contemporary geometry which have their roots in the Age of Enlightenment. {\it Projective geometry} was primarily developed in the 19th century.  Its roots stemmed from works of Pascal from 1639 and from Desargues in 1642, and their contributions were rediscovered by the 19th century geometers%
\footnote{This development of projective geometry was discussed in some detail in our paper \cite{wells2013}.}%
.  {\it Algebraic topology} was hinted at in a letter of Leibniz to Huygens in 1679 \cite{leibniz1679}, which was cited by Euler in his famous paper on the K\"{o}nigberg bridges \cite{euler1735}.  Leibniz and then Euler used the phrase "analysis situs" to describe  a relationship between geometric and algebraic quantities.  This terminology was then used by Poincar\'{e} in his groundbreaking series of papers at the end of the 19th century which created modern algebraic topology%
\footnote{These developments in algebraic topology played an important role in the development of complex geometry in the 19th centurey as we discuss in our paper \cite{wells2015}}%
.  The final paper worth mentioning from this period, and  which is very important for algebraic topology, was Euler's singular paper \cite{euler1752} from 1752 which described for the first time what has become known as the {\it Euler characteristic} for a surface.

We will consider in more detail in this paper  two  other important areas of geometry developed in this time period, namely {\it algebraic geometry} and {\it differential geometry}, both of which had substantive growth in the 18th century.
Section \ref{sec:algebraic-geometry} describes some of the work of Pappus which played an important role in these  geometric developments.  We then turn to Descartes, whose work includes a solution to a problem posed by Pappus, and included the famous coordinate geometry which transformed mathematics in so many ways.  Both Descartes and Fermat were able to successfully classify the algebraic curves of degree two, as described in Section \ref{sec:degree-two}.  Fermat's work on number theory is much more well known, but he contributed significantly to the study of curves in the plane.

Following up on the work of Descartes and Fermat, Newton gave a detailed classification of real-algebraic curves of degree three in the plane.  Some 50 years later Euler  gave a similar classification. Here analysis played an important role in both the work of Newton and Euler, in particular in their use of infinite series in their descriptions of asymptotic behavior. Section \ref{sec:degree-three} gives an brief overview of this initial work of Newton.

 In Section \ref{sec:differential-geometry} we note how the study of transcendental functions led to many geometric objects which were not necessarily defined algebraically (this was pointed out quite explicitly in Euler's {\it Introductio} in 1748 \cite{euler1748}). Earlier, Newton first formulated the notion of curvature of a curve in the plane in terms of calculus, which followed up on ideas of Apollonius who looked at curvature of conic sections  and the more general work of Huygens from the 17th century.  These ideas are all discussed in Section \ref{sec:curvature-plane}.  The curvature of curves in space was initiated by Clairaut in the early 18th century and brought into its final form in the mid-19th century by Frenet and Serret (Section \ref{sec:curvature-space}).  The final topic in this section considers the work of Euler from 1767 who studied the curvature of a surface by analyzing the curvature of the curves arising from intersections of planes normal to the surface with the surface itself.

In the Conclusion, we note that these ideas all relate to a number of 19th and 20th century geometric ideas which involved new elements of abstraction  such as intrinsic differential geometry, projective geometry, and abstract higher dimensional manifolds.

\section{Algebraic Geometry}
\label{sec:algebraic-geometry}
One of the last major figures in Greek mathematics was
Pappus of Alexandria ({\it c.} 290 AD--350 AD).  He published a work entitled {\em Collection}.  Book I and the introduction to Book VII of this major work have been lost, but that which has been preserved gives a good survey of many mathematical discoveries of his predecessors which have been lost, and, in addition, he contributed significantly to solutions of a number of geometric and arithmetic problems (see e.g., Boyer \cite{boyer}, pp. 205-213 for a summary of the important contributions in the {\em Collection}).  We want to single out one particular contribution that has played such an important role in the history of geometry.  This is now referred to as the {\em Problem of Pappus} and was described in the beginning of Book VII of the {\em Collection}.  This problem was treated by Euclid, Apollonius
 and others who preceded Pappus, and, as will be seen later, Descartes.

Apollonius ({\it c.} 262 BC- {\it c.} 190 BC)
 wrote in the preface to his {\em Conics} \cite{apollonius1896} the following (using the translation in Boyer \cite{boyer}, p. 167):
\begin{quote}
The third book contains many remarkable theorems useful for the synthesis of solid loci and determinations of limits;  the most and prettiest of these theorems are new and, when I had discovered them, I observed that Euclid had not worked out the synthesis of the locus with respect to three and four lines, but only a chance portion of it and that not successfully: for it was not possible that the synthesis could have been completed without my additional discoveries.
\end{quote}
Here Apollonius was referring to his discoveries concerning conic sections, which transcended substantially the work of the earlier Greek mathematicians. The "synthesis of the locus with respect to three and four lines" is a special case of what has come to be called the Problem of Pappus which is discussed more explicitly in the following paragraphs. 

We first formulate the original version, which had been solved by Apollonius, and then later its generalizations. The Problem of Pappas: given three or four lines in the Euclidean plane, find the locus of points such that the the square of the distance to one line (in the three-line case) is proportional to the product of the distances to the remaining two lines.  In the case of four lines, one asks for the locus of points with the property that the product of the distances to two of the lines is proportional to the product of the distances to the remaining two lines%
\footnote{This set of problems is similar in spirit to the characterization of a circle being the locus of all points (in the plane) whose distance to a given point is constant or an ellipse being the set of all points such that the sum of the distances to two distinct points is constant.}%
. In all cases the distance to a given line is measured at a given angle to the given line (thus the given data is the set of lines and the set of angles and the proportionality factor). Apollonius  shows that the resulting curves are indeed given by conic sections, which, of course, is the primary topic of his book. He implies in the quote above that Euclid did not have the detailed results needed concerning conic sections in order to solve this problem, which is very likely the case. 

Using the language and ideas of analytic geometry, one can easily verify Apollonius' result (see Boyer \cite{boyer}, pp. 167-168). If, in the case of three lines given by the equations
\bean
A_1 x + B_1y + C_1  & = & 0,\\
A_2 x + B_2 y+ C_2 &=& 0,\\
A_3 x + B_3 y + C_3 & = & 0,
\eean
and if the angles used for measuring distance are given by \(\theta_1\),  \(\theta_2\), and \(\theta_3\), then the locus is given as the set of points \((x,y)\) satisfying:
\[
\frac{(A_1x+B_1y+C_1)^2}{(A_1^2+B_1^2)\sin^2\theta_1}=K\frac{(A_2x+B_2y+C_2)}{\sqrt{A_2^2+B_2^2}\sin \theta_2}\cdot \frac{(A_3 x + B_3 y + C_3)}{\sqrt{A_3^2+B_3^2}\sin \theta_3}.
\]
Since the locus is the solution of a quadratic equation in the plane, it follows that it is a conic section, which is what Apollonius
 had discovered using his methodology. 

The general problem of Pappus is to be given an arbitrary number of lines and angles and to ask the same question.  Here is a quote from Pappus concerning the more general problem.  First he notes (following Boyer \cite{boyer}, p. 209) that for six lines, the locus can be considered that a solid is in fixed ratio to another solid (here ``solid'' refers to products of three lengths, i.e., homogeneous polynomial terms of degree three using modern language) However, higher degree terms were a mystery to him as (quoting Pappus)
\begin{quote}
there is not anything contained by more than three dimensions
\end{quote}
and, he continued
\begin{quote}
men a little before our time have allowed themselves to interpret such things, signifying nothing at all comprehensible, speaking of the product of the content of such and such lines by the square of this or the content of those. These things might however be stated and shown generally by means of compounded proportions.
\end{quote}
Pappus did not study the higher degree case (higher than six), but he did make the important observation that the the loci were curves in the plane.  As Boyer observes, Pappus was a geometer and Diophantus, a contemporary, was an algebraist (who did consider higher powers and who had the notation to tackle the higher degree problems), but it required a mathematician who was familiar with both algebra and geometry to make the next step, and that turned out to be Descartes, some 1300 years later.

\subsection{Algebraic Curves of Degree Two: Descartes and Fermat}
\label{sec:degree-two}
Ren\'{e}  Descartes (1596-1650) published a slim volume in 1637 entitled {\em La G\'eom\'etrie} \cite{descartes}, which initially was an appendix to a longer work in philosophy but was also published independently. Descartes's work turned out to be revolutionary, and  when the next generation of mathematicians began to write general texts concerning what today is called analytic geometry  the impact of his work spread throughout the mathematical world of Europe and became fully developed in the 18th century.  Until Descartes, and actually long after as well, Euclid's {\em Elements} were definitive on almost all things concerning geometry. Descartes' new view of geometry was very important, but Euclid's ideas were still very valid.  Only with non-Euclidean geometry in the 19th century was Euclid challenged in a fundamental way. A very brief but succinct survey of Descartes' {\em G\'eom\'etrie} is given by Serfati \cite{serfati1}.

Descartes is most well known to mathematicians for having discovered {\em analytic geometry} or probably more appropriately named {\em coordinate geometry}, which has been taught in twentieth century high schools and on up to the present day around the world.  The label ``analytic geometry" as applied to Descartes is slightly a misnomer. What Descartes showed was how common problems of geometry as described by Greek geometers could be described by using algebraic equations and conversely.  The most important historical example of this is Descartes' theorem that solutions of algebraic equations of degree two in two variables corresponds precisely to the conic sections studied by Apollonius
 and others.  What was much more important, aside from the  new coordinate system point of view, was that he {\em defined} geometric objects to be solutions of algebraic equations. In particular, he considered algebraic equations in two variables of arbitrary degree, which became the nucleus of  {\em algebraic geometry}, and which, in its modern form as developed in the 19th and 20th century, is definitely {\em not} taught in high schools around the world.

\begin{figure} 
\centerline{\begin{tabular}{cc}\includegraphics[width=6cm]{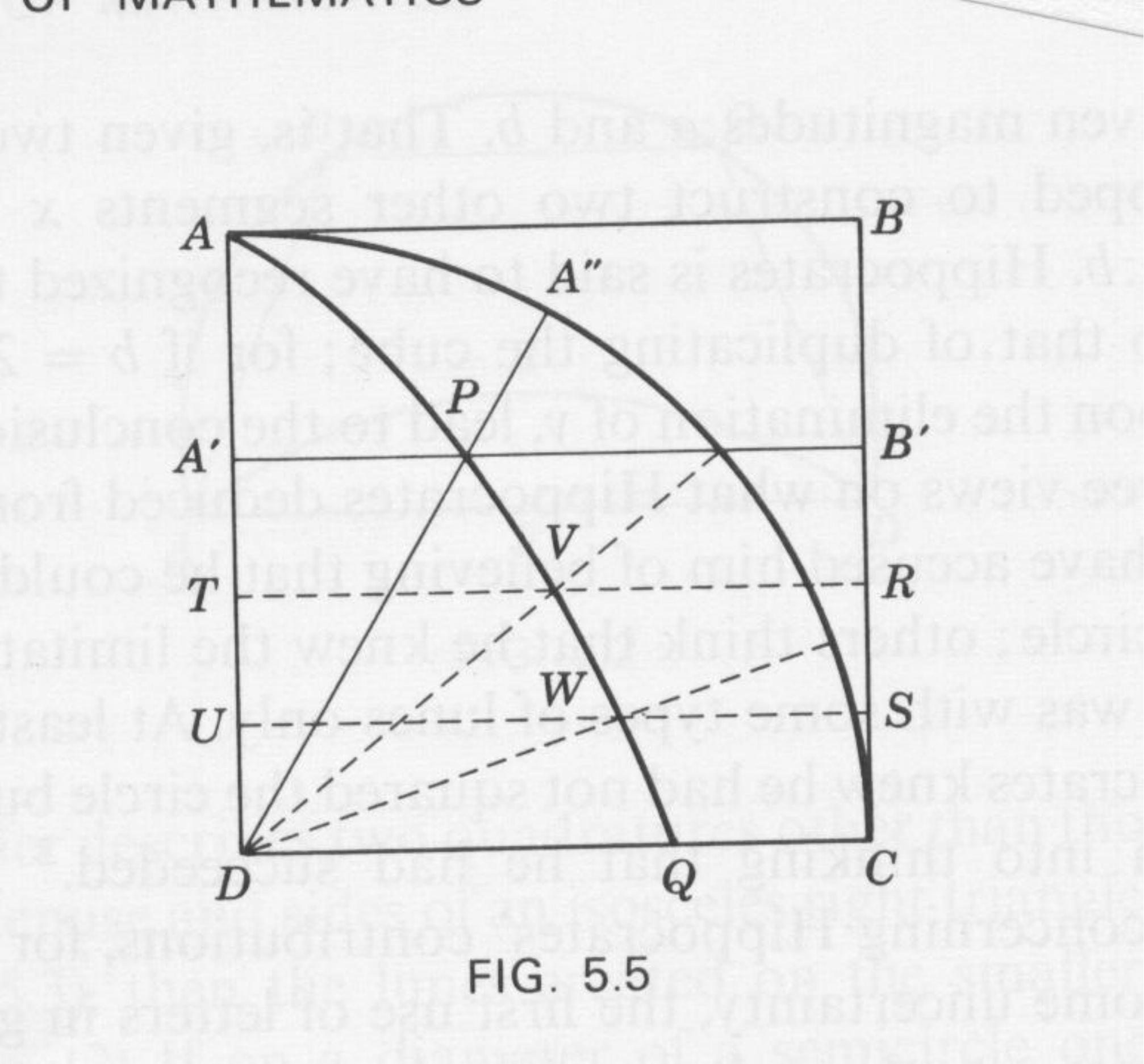}&\includegraphics[width=6cm]{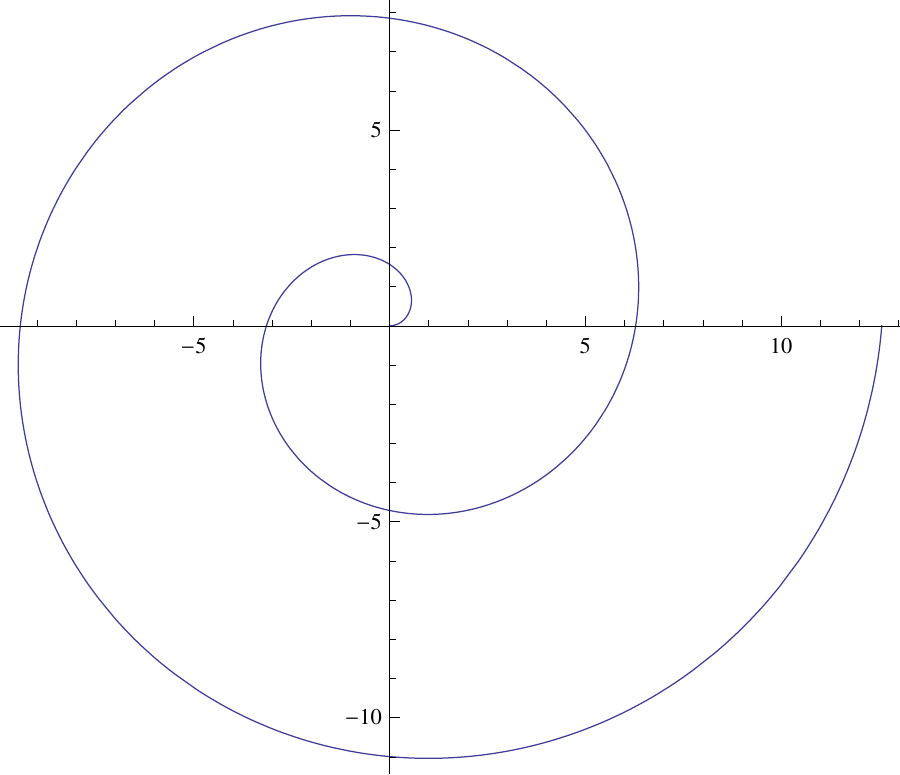}\end{tabular}}
\caption{The quadratrix (trisectrix) curve of Hippias of Elis and the spiral of Archimedes. In polar coordinates, \(r =\frac{2a\theta}{\pi \sin \theta}, r=a\theta\), respectively.}
\label{fig:quadratrix}
\end{figure}

However, the way analytic geometry is taught today involves not only algebraic functions, but also the standard transcendental functions, $\sin x$, $e^x$, etc., as well, and students learn about the curves that these functions can represent in the plane as well.  This is something that Descartes absolutely rejected.  He had learned from Greek authors that there were three types of curves studied by mathematicians: {\em plane curves}, that is curves that could be described with a straight edge and compass in the plane; {\em solid curves}, that is curves that could be described in three space by intersections of simple surfaces with a plane, the simplest being the full family of conic sections, i.e., intersections of a cone with a plane; and the third category was {\em linear curves}, i.e., everything else. This last category included the quadratrix (often called the trisectrix) , the spiral of Archimedes, and other transcendental curves (to use modern language). These are illustrated in Figure \ref{fig:quadratrix}.
  Of course, this terminology (which one finds in Euclid, Apollonius
 and Pappus) would be totally confusing today.  Descartes pleas for his reader to have only two categories of curves: {\em geometric} and {\em mechanical}.  His definition of geometric was simply anything described by algebra, i.e., {\em algebraic} curves (and surfaces, which he did not really address in his book), and {\it mechanical} being all others.  Hence he extended the Greeks planar and solid curves to include curves of arbitrary degree.  On the other hand, he excluded from the study of geometry the mechanical curves which include the quadratrix and spiral.  The conchoid is an example is an algebraic curve of higher degree
\[(x-a)^2(x^2-y^2)=b^2x^2\]
and this curve was used for cube duplication and angle trisection problems. The squaring of the circle required trigonometric functions (for instance, using the quadratrix). Today one uses the terminology of Leibniz: {\it algebraic curves} and {\it transcendental curves}.

The main reason that Descartes made this distinction was that he thought that one could always calculate solutions to algebraic equations but not solutions of transcendental equations. He was familiar with the solutions for the third and fourth degree formulas of Cardano, and it seems that he assumed that such formulas were true in general.  Descartes spent almost all of Book III (the third and final chapter of {\em La G\'eom\'trie}) discussing the explicit solution of equations and in particular he developed a theory of roots of polynomial equations in one variable.  But in any event, if he had a polynomial equation of degree $d$ of the form $P(x,y)=0$, he argued that if you fixed a particular value of the variable $y$ then you obtained a polynomial in one variable which could always be solved (he believed the fundamental theorem of algebra was indeed true, but he also believed roots could be found simply by extracting roots, which Abel and Galois later showed to be false). He then argued that one could not compute the values of the function whose graph is the quadratrix curve for arbitrary values of the variable in question, but only for certain special values.  Fundamentally, this is equivalent to the fact that one can compute by elementary geometry special values of the trigonometric functions, but not all of them.   This was beyond the scope of Descartes at the time, and for this reason he rejected the study of transcendental functions as an object of study in geometry. In the eighteenth century this would change radically.  Of course Descartes excluded solutions of equations which involved complex numbers with nonzero imaginary part, and said simply there were no solutions in those cases.  It would be two more centuries before mathematicians became comfortable with complex numbers (we discuss this in some detail in \cite{wells2015}).

Descartes had strong opinions on what was true (or worthy of study) and what was not.  One more example, which is of historic importance, concerns arc length.  In the development of trigonometric functions, arc length is a critical ingredient (the relation between the length of an arc to the length of the chord subtended by it is how trigonometry was originally introduced, and the sine and cosine functions are simply modern variations of this). Descartes was convinced that for no algebraic curve (e.g., the circle) could one ever find precisely the length of the arc in terms of the length of the chord (this clearly relates to the difficulty of computing $\pi$). As he put it on p. 32 of {\em La Geometrie} 
\begin{quote}
...car encore qu'on n'y puisse recevoir aucunes lignes qui semblent \`{a} des cordes, 
c'est-\`{a}-dire qui deviennent 
tant\^{o}t droites  et tant\^{o}t courbes, \`{a} cause que la proportion
qui est entre les droites et les courbes n'\`{e}tant pas connue, et m\^{e}me, je crois, ne
le pouvant \^{e}tre par les hommes, on ne pourroit rien conclure de l\`{a} qui f\^{u}t exact
et assur\'{e}.%
\footnote{
... because one should not be able to consider lines (or curves) that are like strings, in that they re sometimes straight and sometimes curved, since the ratios between straight and curved lines are not known, and I believe cannot be discovered by human minds, and therefore no conclusion based upon such ratios can be accepted as rigorous and exact.  }
\end{quote}

What Descartes did do, and it seems to have been a major part of the inspiration for writing his book, was to give a new and extensive solution to the Problem of Pappus (both Books I and II are devoted to this topic, among other things). He showed that for an arbitrary number of straight lines and associated angles in the plane the locus of points such that the products of the distances to half the set of lines was proportional to the products of the distances to the other half of the set of lines, all distances being measured at the given angles,  is an algebraic curve in the plane%
\footnote{For the case of an odd number of lines, one takes the distance to one of the lines twice in this proportionality.}%
.   He computes that, for the classical case of three or four lines (the original problem solved by Pappus), the curve is an algebraic curve of degree 2, for the case of 5--9 lines the curve is an algebraic curve of degree 4, and for 10--13 the curve is of degree 6, etc.  He refers to the curves of degree 2, 4, and 6, etc., as curves of {\em genre} 1,2, and 3, etc. Descartes considered various special cases where odd degree polynomials could appear, but he lumped them in his genre classification with the even degree cases.

But, in addition to showing that the solution to the problem was algebraic curves, he showed that for the classical case where he got algebraic curves of degree two, that all of these curves were indeed conic sections. In his proof he showed how each of the polynomials of degree two arising in this context could be put in a normal form by a suitable change of coordinates, and then he was able to use Apollonius'
\index{Apollonius ({\it ca.} 262 BC- {\it ca.} 190 BC)}%
 characterization of the conic sections in terms of suitable coordinates to determine that the solutions were indeed conic sections. The major distinction between Apollonius
\index{Apollonius ({\it ca.} 262 BC- {\it ca.} 190 BC)}%
 and Descartes in this context was that Apollonius
\index{Apollonius ({\it ca.} 262 BC- {\it ca.} 190 BC)}%
 started with a given conic section and produced coordinates which helped describe it, while Descartes started with the coordinate system and the equation and was able to put it in canonical form and identify it in a suitable manner.

Descartes recognized that different equations could describe the same geometric curve, and he pointed out the need to find the ``simplest" algebraic function that could represent a given curve.  This, of course, was the key question for the classification of algebraic curves (and later higher dimensional manifolds in a variety of categories, to use a small pun) which has been a consistent and important theme in the following centuries.  Descartes classified the algebraic curves of degree two, and Newton followed up with his major work on the classification of algebraic curves of degree 3 \cite{linearum} which 	 will be discussed in more detail in the next section.

Pierre de Fermat (1601-1665) played an important and less recognized role in this development of geometric ideas. First, in his only published paper in his lifetime%
\footnote{See  \cite{mahoney}, p. 267 for this reference {\it De linearum curvarum cum lineis rectis comparatione dissertatio geometrica} (Geometrical dissertation on the comparison of curved lines with straight lines), which appears as an appendix in a book by Antoine de Lalouv\`{e}re from 1660).}%
, he showed that one could explicitly compute the arc length of a specific algebraic curve, which , as noted above, Descartes claimed was impossible. More precisely, he showed that for the algebraic curve $y^2=x^3$, the semicubical parabola, the arc length could be explicitly computed.  Namely, if one takes the positive branch of this curve $y=x^{\frac{3}{2}}$ on the interval $[0,b]$, then one can verify that the arc length integral is
\[
\int_0^b\sqrt{1 + [y'(x)]^2}dx = \frac{4}{9}[(\frac{4}{9})^2(b-1)^{\frac{3}{2}}]
\]
as any calculus student today can do (and is often asked to do!).  However, Fermat, who played a major role in the development of differential and integral calculus, was not aware of the fundamental theorem of calculus or of this arc length formula, which makes the calculation somewhat more difficult!

\begin{figure} 
\vspace{6pt}
\centerline{
	\includegraphics[width=13cm]{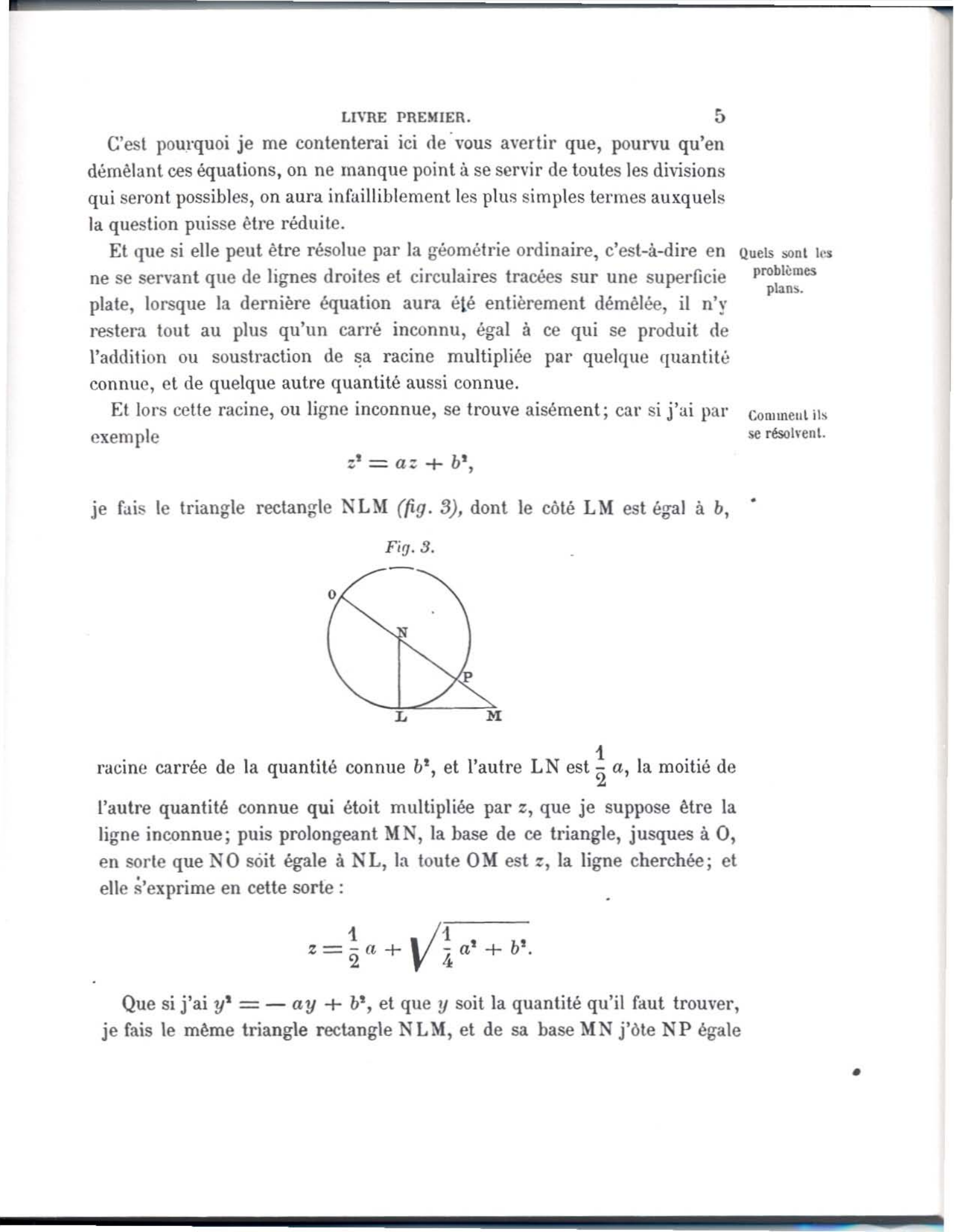}}
 \caption{Page 5 of Descartes' {\em La G\'eom\'etrie}
}
\label{fig:descartes-p5}
\end{figure}

Moreover, Fermat also proved, independently of Descartes, that algebraic curves of degree two were conic sections (see his biography which discusses this among many other things \cite{mahoney} (see also \cite{resnikoff-wells1984}).
A major difference between these two historic figures on this particular point was that Fermat expressed his work in the classical language of Euclidean geometry $\overline{OA}$, $\overline{OB}$, etc. representing the lengths of the line segments $O$ to $A$ or $O$ to $B$ in the Euclidean plane. Descartes used this notation as well, but he adroitly introduced the variables $x$, $y$, and $z$ to denote unknown quantities (lengths of segments in the problems he was considering), and symbols $a$, $b$, $c$, etc., from the first part of the alphabet to represent known quantities in a given computation, a practice that has been followed ever since.  In this way he reduced geometric problems to algebraic problems.  He was very concerned that his new way of looking at things should be well connected with classical Greek geometry.  As an example in Book I of {\em La G\'eom\'etrie} he goes to great pains to show that solutions of a quadratic equation such as
\[
z^2 = az+b^2
\]
can be constructed by straight edge and compass.  Figure \ref{fig:descartes-p5}
illustrates his construction.

We want to mention one final historical note that is not that well known to the mathematical public (or public in general).  One of the most important innovations in Descartes' {\em La  G\'eom\'etrie}  was the invention of {\it exponential notation}. He used $x^3$, $x^4$, etc. freely throughout the book.  In the 1886 edition that we have been quoting from, the editors point out two modernizations that they introduced to make the book more readable for the nineteenth century reader. The first was the use of the $``="$ sign instead of the symbol $\infty$ that Descartes used for equality, and the second, surprisingly, was the use of $x^2$ instead of Descartes' preferred notation $xx$.  This certainly seems strange to a modern reader, as he used the higher power exponential notation with no hesitation.  Before this innovation of Descartes mathematicians used {\em different symbols} for different powers of the unknown variable $x$, which would make a formula like the law of exponents
\[
x^{m+n} = (x^m)(x^n)
\]
somewhat difficult to formulate (see, e.g., \cite{resnikoff-wells1984} and the references therein).

\subsection{Algebraic Curves of Degree Three: Newton and Euler}
\label{sec:degree-three}
\begin{figure} 
\vspace{6pt}
\centerline{
	\includegraphics[width=12cm]{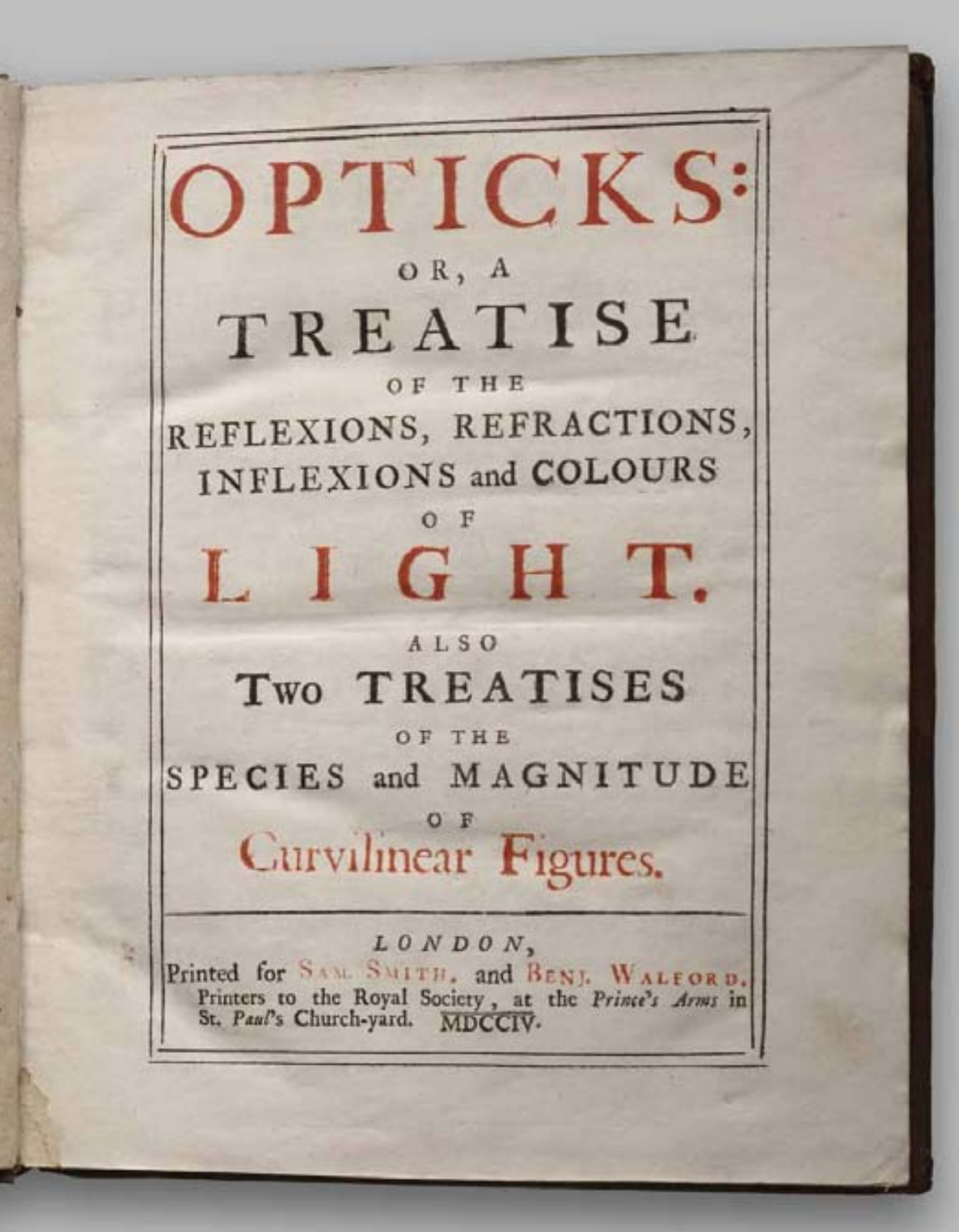}}
 \caption{Title Page of Newton's {\em Opticks} from 1704
}
\label{fig:NewtonOpticsTitlePage}
\end{figure}

The developments of differential and integral calculus and other ideas in analysis (e.g., the theory of infinite series) in the late seventeenth century by Isaac Newton (1642--1746) and Gottfried Wilhelm Leibniz (1646--1716) and their successors has been one of the most important developments in all of mathematics.  It has been well documented (see, e.g. Boyer \cite{boyer}, Kline \cite{kline1972}, and other such references).  We won't try to give any historical background on this important topic, as we want to concentrate on the interaction between analysis with the developments in geometry as it evolved in this Age of Enlightenment. Newton published in 1704 for the first time two pivotal works on geometry and calculus which were quite independent of each other and which were chapters in a larger book of optics in \cite{opticks} . Figure \ref{fig:NewtonOpticsTitlePage}
shows the title page of this very important work; here the two mathematical treatises are labeled "Also Two Treatises of the Species and Magnitude of Curvilinear Figures. The first of these mathematical treatises deals with the classification of algebraic curves of degree three \cite{linearum}, and the second concerns itself with calculus and measuring the area under a curse \cite{quadratura}. Both of these mathematical chapters are written in Latin, whereas the main part of the book is written in English, and deals with physics, and optics in particular. Note that Newton's {\it Principia Mathematica} was published earlier in 1687 and the definitive third edition was published in 1726 (see the annotated copy of this third edition of Newton's most important work \cite{principia} edited by Khoyr\'e and Cohen, which was published in 1972). This very important book uses the ideas of calculus in a fundamental way, even if the language used expresses the calculus ideas in somewhat cumbersome terms of limiting processes of Euclidean geometric objects. His work on algebraic curves is, however, not a part of Principia Mathematica.

Newton gave a quite precise classification of algebraic curves of degree three. This is a direct generalization of the case of curves of degree two, the conic sections. The basic tool used was to analyze the highest order homogeneous terms of degree three and their possible factorizations over the real numbers. This led to various types of branches that are unbounded (and to various types of asymptotes)  that can arise, and they become an important part of the classification. For instance, for curves of degree two, one can see that  an ellipse has no infinite branches, a parabola has two infinite branches, and a hyperbola has four infinite branches including two straight lines which are asymptotes. This behavior at infinity completely distinguishes these three classes of curves.

\begin{figure} 
\vspace{6pt}
\centerline{
	\includegraphics[width=12cm]{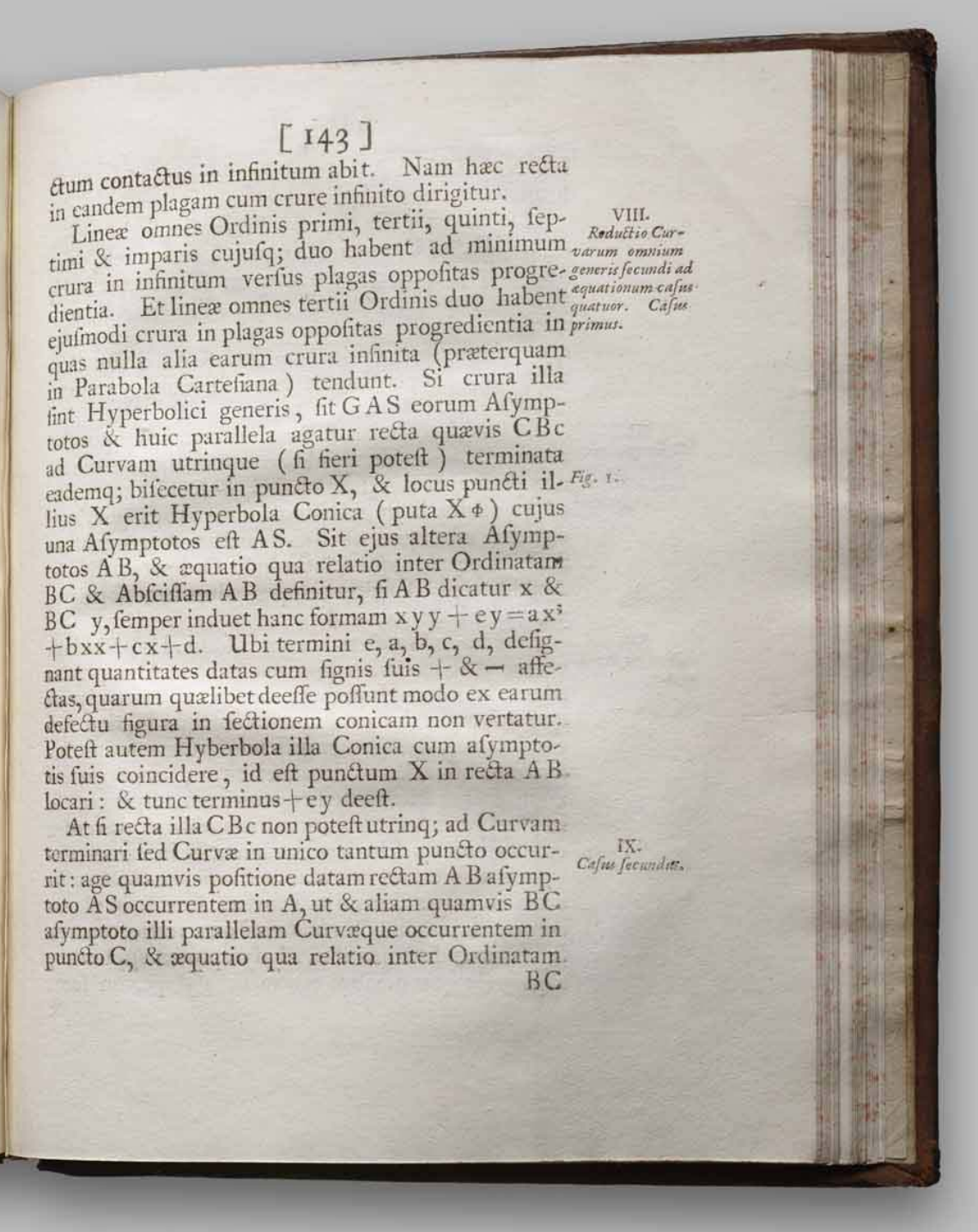}}
 \caption{Page 143 of Newton's {\em Linearum} from 1704
}
\label{fig:newton_p143}
\end{figure}

\begin{figure} 
\vspace{6pt}
\centerline{
	\includegraphics[width=12cm]{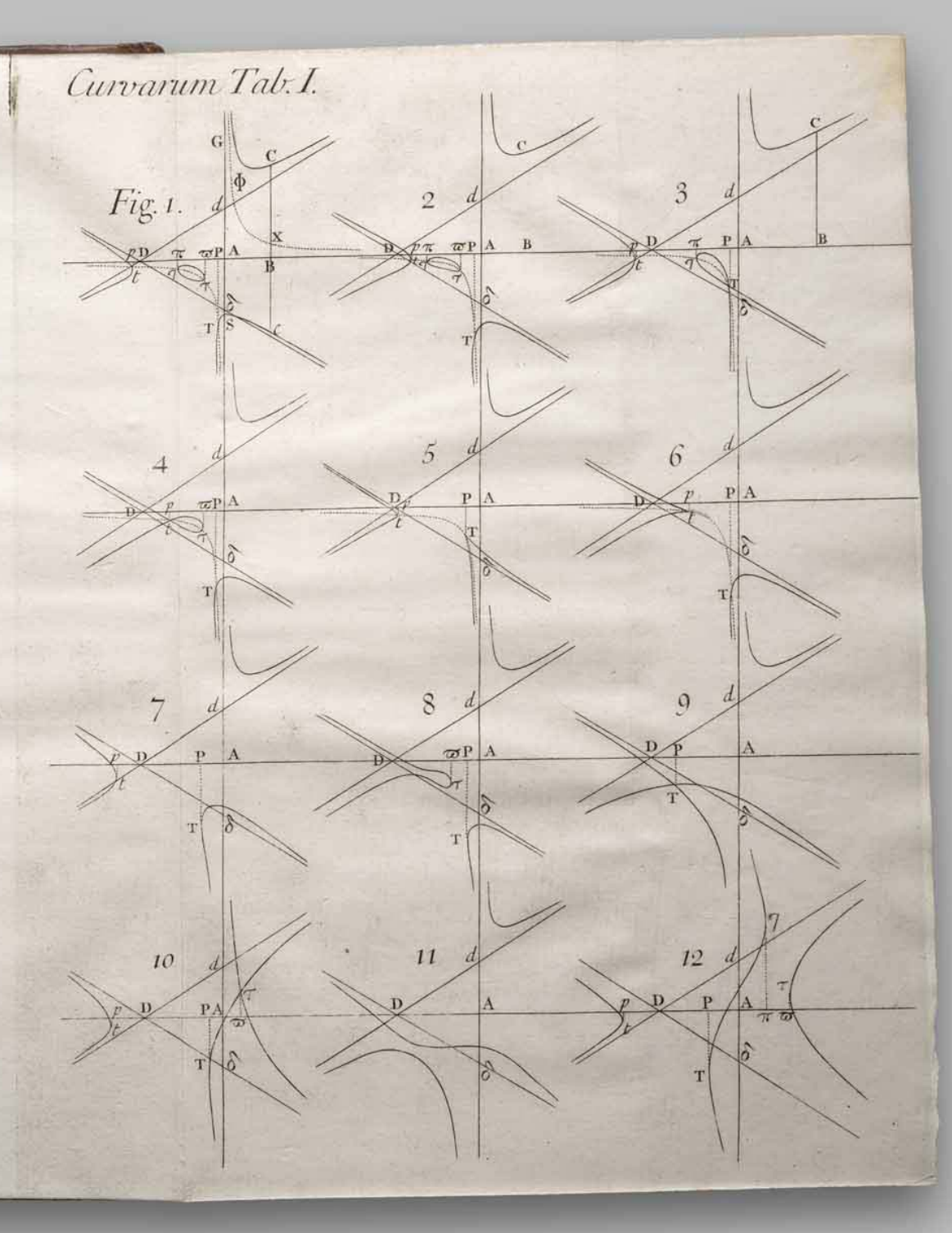}}
 \caption{First page of Newton's figures from {\em Linearum} including Fig. 1 referenced  on p. 143 of the same book}
\label{fig:newton_fig1}
\end{figure}

Newton's analysis yielded the classification of 78 different types of curves. In fact he only described 72, having missed six types. He classified them algebraically, and then provided beautiful drawings of the typical curve of the specified classification type. In Figures \ref{fig:newton_p143} and \ref{fig:newton_fig1}
we see some examples of the classification and the corresponding drawings. His analysis of how he arrived at the classification is very terse indeed. In fact, there is very little description of his analysis.  The monograph is not much more than a simple listing of his findings. Proofs of his results were published later by mathematicians who analyzed and generalized his results. He indicates also that one can carry out such an analysis for curves of higher degree by the same method of analysis

Many mathematicians over the centuries have, of course, analyzed Newton's results in great detail and one can find the classification in a variety of monographs over the centuries since (see for example Walker \cite{walker} for a modern treatment of real algebraic curves). However, one of the most beautiful and thorough analyses of Newton's work from the classical literature is that of Euler in \cite{euler1748}. This very influential book, {\it Introductio in Analysis Infinitorum}, was first published in Latin in 1748 and then reappeared later in many additions and in various languages.  For instance one can buy the book today in English on Amazon.com, and there is a very nice French edition from the late 18th century published at the time of Napoleon \cite{euler1796}.

What we have described above very briefly is Newton's classification of algebraic curves of degree three from the beginning of the 18th century.  Such developments carried on into modern times, in particular with similar classifications of Riemann surfaces (considered as algebraic curfaces complex dimension one  of various degrees), and complex manifolds of higher complex dimension. In general the study of complex manifolds or complex algebraic varieties of various dimensions has turned out to be simpler than the study of real algebraic manifolds and varieties over the real numbers due to the closure of the field of complex numbers. But it all started in the real algebraic setting, as that's what the mathematical community was familiar with in the early 18th century. They knew about complex numbers, but they were not yet familiar with complex geometry.

\section{Differential Geometry}
\label{sec:differential-geometry}
\label{sec:differential-geometry}
A second major development in geometry in the 18th century was the study of curves and surfaces in \(\BR^2\) and \(\BR^3\) defined by not necessarily algebraic functions. These included two not quite independent developments that took place more or less simultaneously. The first was the development of the now standard elementary transcendental functions: the trigonometric, exponential, and logarithmic functions. In Euler's textbook from 1748 \cite{euler1748} these functions and their algebraic and analytic properties (e.g., \(\frac{d}{dx}\sin (x) =\cos (x), \sin (x+y) = \sin (x) \cos (y) + \cos (x) \sin (y)\), etc.) were fully developed and correspond to what one learns in contemporary precalculus and calculus courses in high school today. The second development involved the solution of differential equations (primarily ordinary differential equations) which provided a large variety of functions for analysis and geometrical representation This led to a large class of special functions that went by the names of the mathematicians who created and developed them: Hermite, Legendre, Bessel, Euler's Gamma function and many others.  These various functions were tabulated for computational use and their various algebraic and analytical properties were developed, similar to those properties illustrated above for trigonometric functions. Over the course of time these mathematical tools became very important for the applications of mathematics to the worlds of chemistry, physics, biology and other areas of scientific understanding. These tools preceded by one or two centuries contemporary methods for scientific analysis available through the use of computers and simulation tools involving modern numerical analysis, which led to the role of special functions being not quite as important as they once were. 

In the latter half of the eighteenth century differential geometry of curves and surfaces began to develop and flourish. First we consider the development of what became known as planar curves and space curves (i.e., smooth curves in \(\BR^2\) and \(\BR^3\)).  {\it Differential geometry} was named as a concept by Bianchi in 1894 (as noted by Kline  \cite{kline1972} on p. 554). This naming of the discipline came long after the most significant developments in the field. It came to mean precisely manifolds equipped with a Riemannian (or more general) metric, or more generally a connection, and where the concepts of curvature played a central role.  Indeed the interaction of differential analysis (i.e. calculus, differential equations, all aspects of analysis involving infinite processes) with geometry is much older and indeed broader than the more precise notion of differential geometry as it is employed today. For instance the notion of differential topology, which developed in the mid-twentieth century certainly involves manifolds and analysis, but doesn't formerly use the notion of a differential-geometric metric as in differential geometry per se. Archimedes knew how to compute areas by the method of exhaustion, and Fermat understood both differentiation of functions (finding maxima and minima and tangents), and how to compute area under some curves, but he did not know the fundamental theorem of calculus (see \cite{resnikoff-wells1984} for a discussion of these issues).  
All of these are indeed an interaction of analysis with geometry, and are parts of the foundation of what became differential geometry two centuries later. 

\subsection{Curvature of curves in the plane}
\label{sec:curvature-plane}
The first important task in differential geometry was to be able to  efficiently compute the tangent line to a given curve at a given point and, as any beginning student of calculus knows, this is one of the first applications of the notion of a derivative. A deeper question, that we explore in greater detail in this section is: what is curvature? More precisely, what is the curvature of a curve in a plane or in three-dimensional space? What is the curvature of a surface in three-dimensional space? Finally, what is curvature on an abstract two-dimensional or higher dimensional manifold?  This last question  is a key part of the geometric developments in the 19th century, and won't be discussed in this paper%
\footnote{We explored this question in some detail in our paper \cite{wells2013}.}%
.

Consider first the simple case of a curve in the plane defined by the graph of a function as in Figure \ref{fig:curvature1}, 
\begin{figure} 
\vspace{6pt}
\centerline{
	\includegraphics[width=12cm]{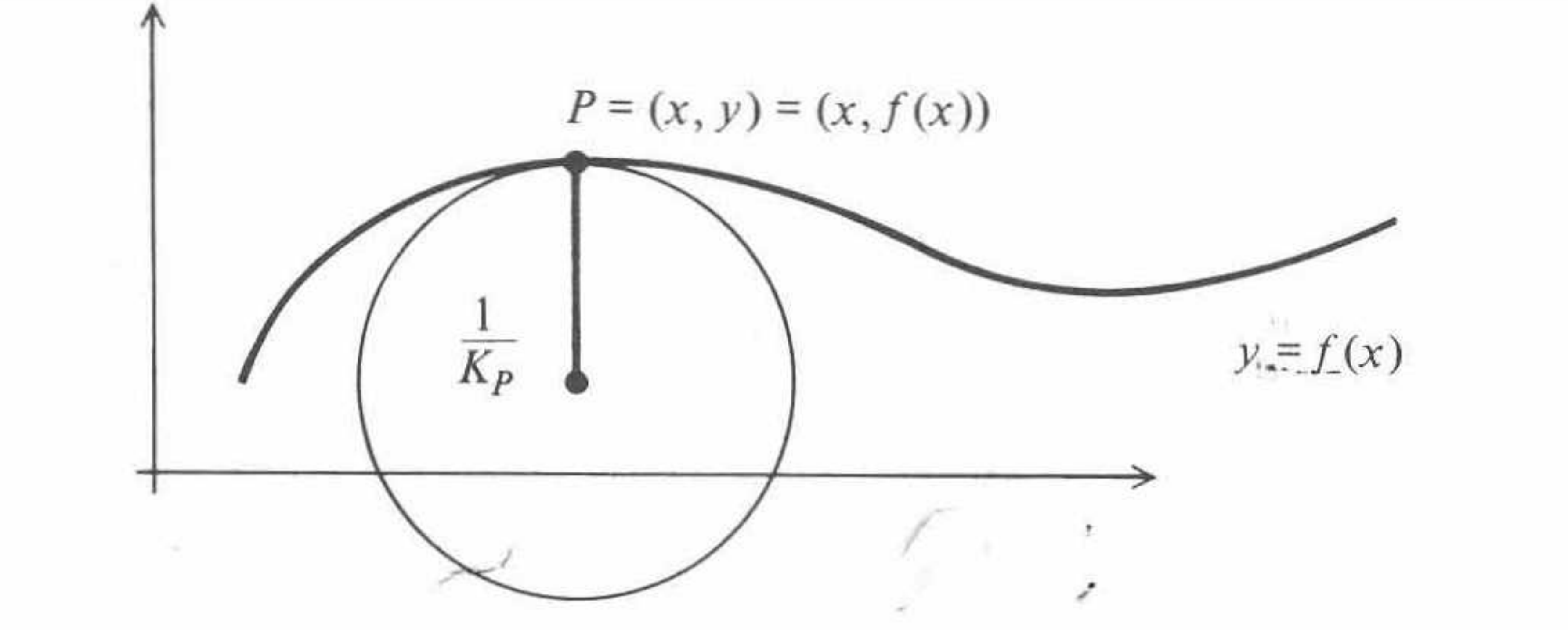}}
 \caption{Radius of curvature of a curve at a point}
 \label{fig:Curvature1}
\end{figure}
then one learns in calculus that the curvature of the curve at \(P=(x,y)\) is given by 
\be
\label{eqn:plane-curvature}
K_P=\pm \frac{f''(x)}{[1+(f'(x)^2]^\frac{3}{2}},
\ee
where the sign is chosen to be positive if the normal vector to the curve at \(P\) intersects the approximating circle and is negative otherwise.  In the illustration in Figure \ref{fig:curvature1}, the normal vector to the curve at \(P\)  using the usual orientation would be pointing upwards in the figure, away from the approximating circle, whose radius is \(1/|K_P|\), and hence in this case the curvature would be negative.

\begin{figure} 
\vspace{6pt}
\centerline{
	\includegraphics[width=12cm]{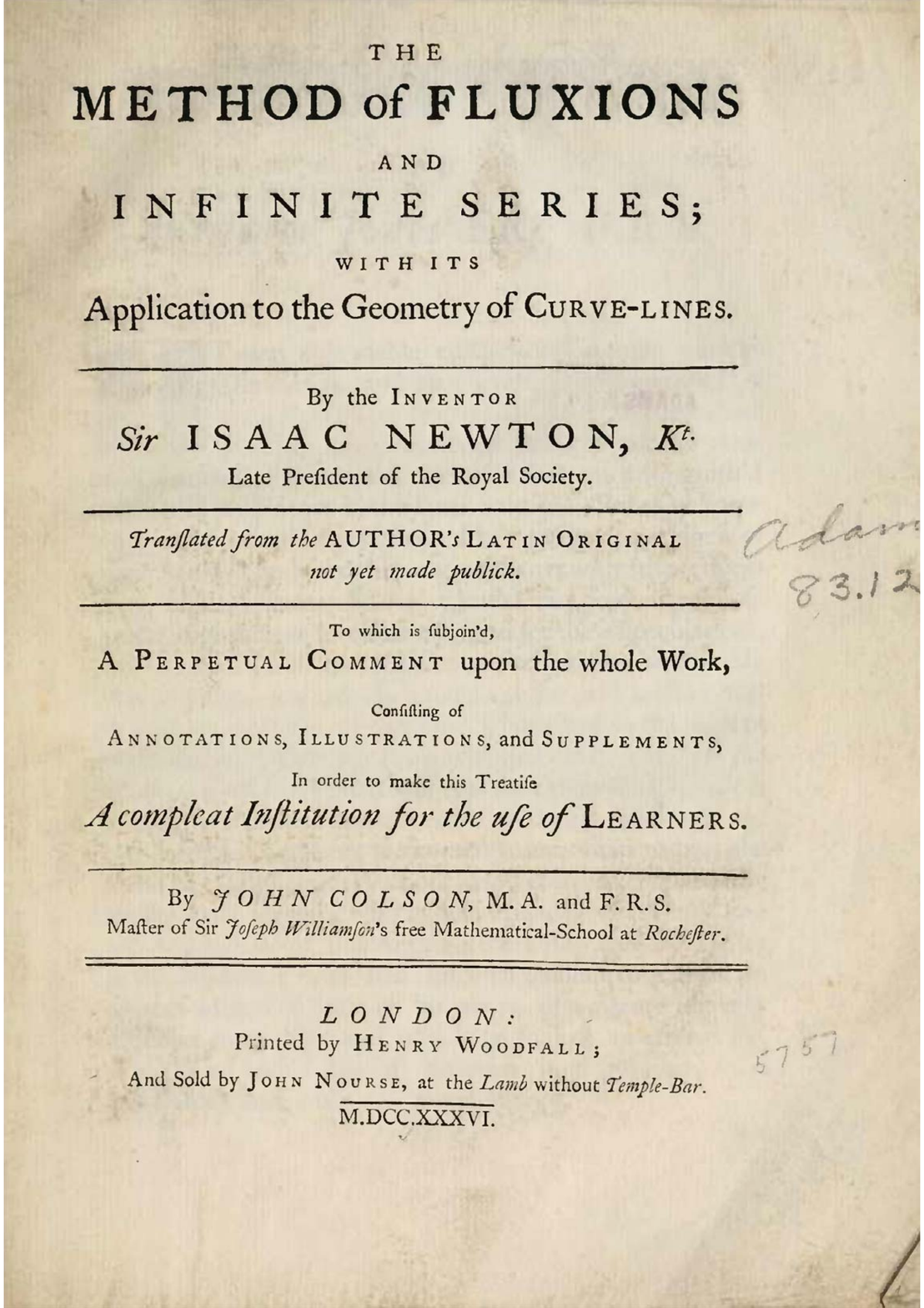}}
 \caption{Title page of Newton's 1736 Monograph on Fluxions}
 \label{fig:newton1736coverpage}
\end{figure}

\begin{figure} 
\vspace{6pt}
\centerline{
	\includegraphics[width=12cm]{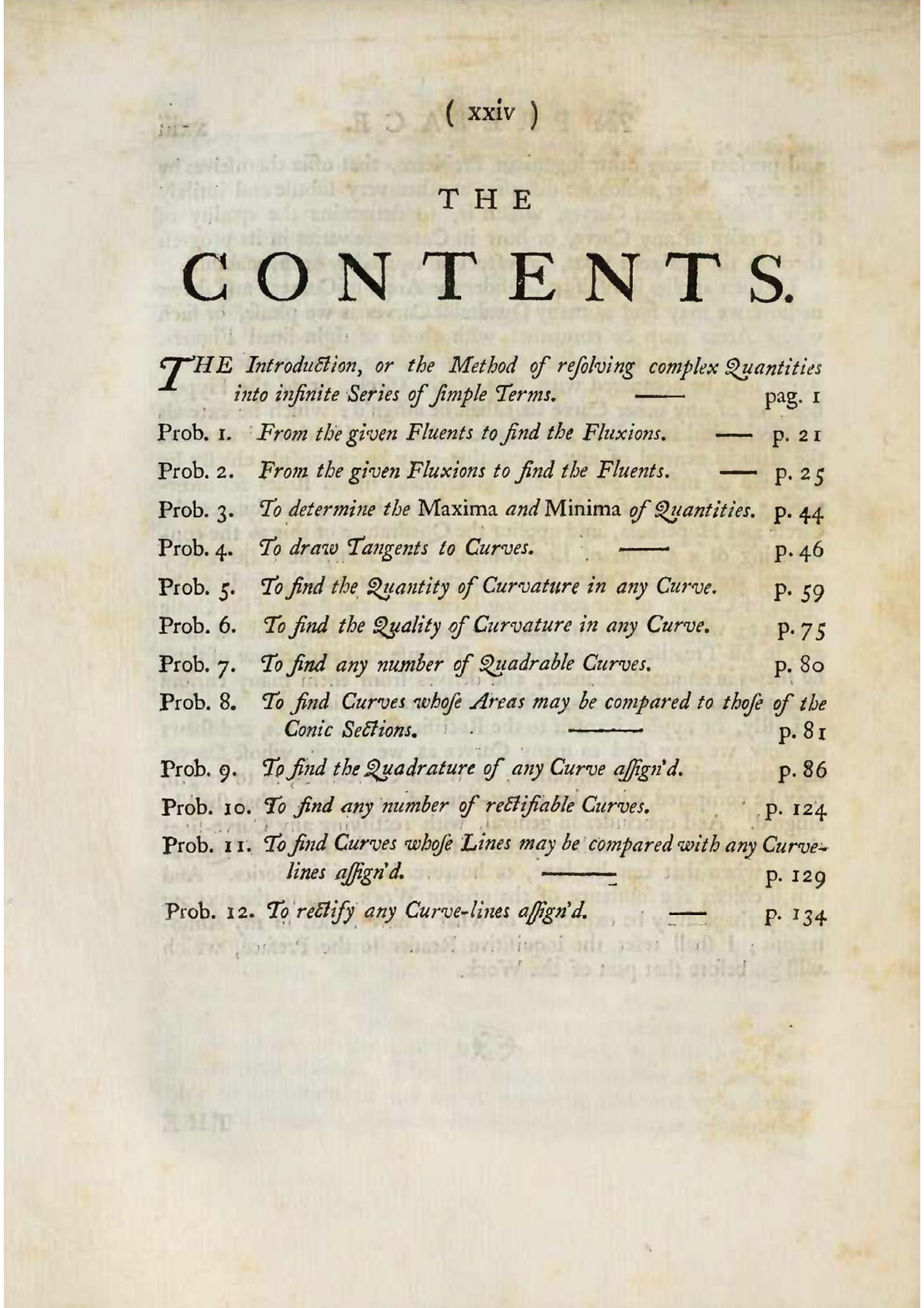}}
 \caption{Table of Contents of Newton's 1736 Monograph on Fluxions}
 \label{fig:newton1736TOC}
\end{figure}

This formula is given for the first time in Newton's monograph of 1736 \cite{newton1736}, which was published as an English translation of his original Latin manuscript from 1671 which was never published, but was privately circulated among some of Newton's colleagues.  This monograph, published in 1736 after Newton's earlier death was part of the basis for the controversy between Newton and Leibniz on who had first invented (or discovered) calculus. Figure \ref{fig:newton1736coverpage}
 shows the cover page of this singular monograph and Figure \ref{fig:newton1736TOC}
 shows the table of contents, where the curvature of a curve stands out so very distinctly as an object of study.

\begin{figure} 
\vspace{6pt}
\centerline{
	\includegraphics[width=12cm]{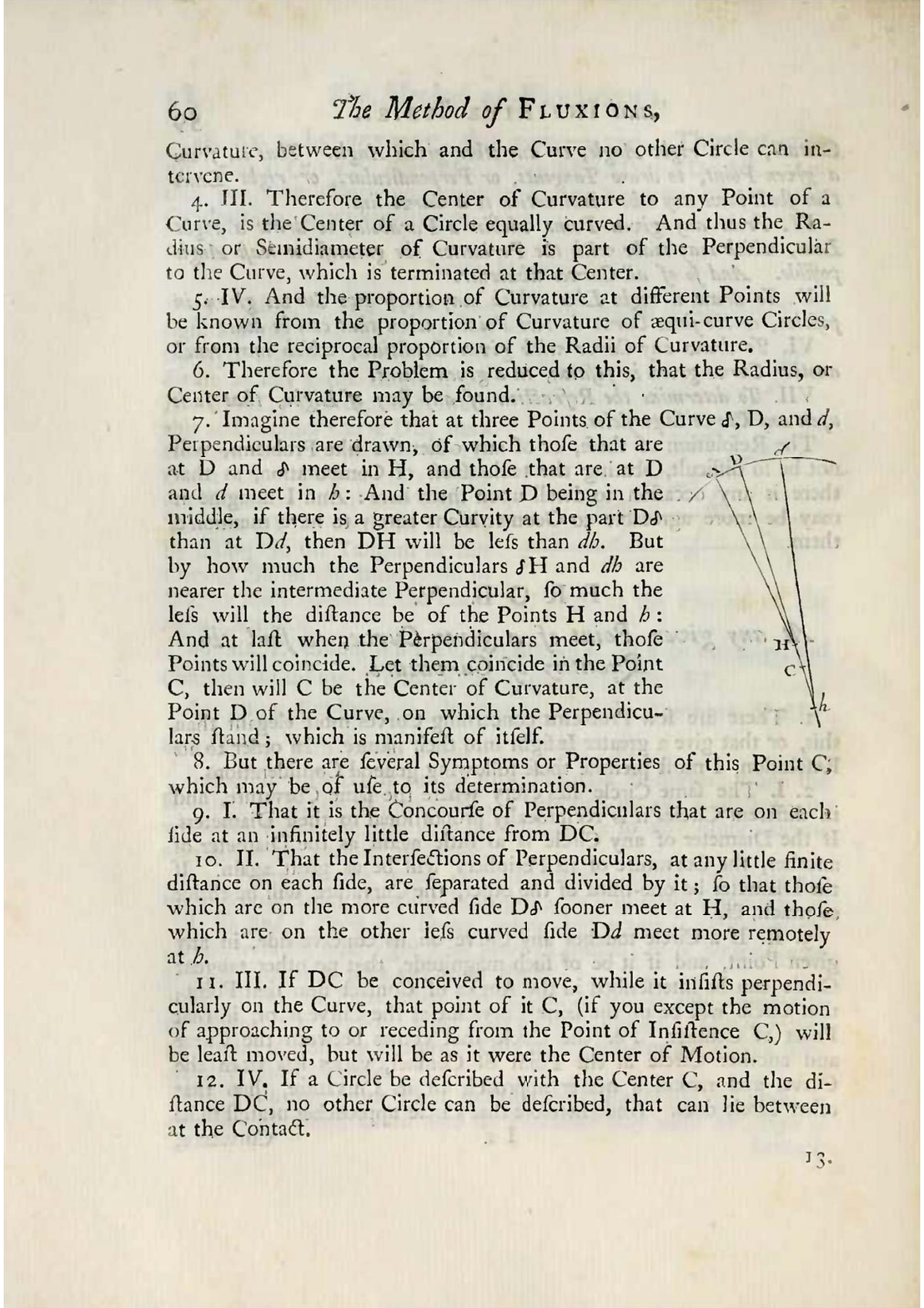}}
\caption{Newton's center of curvature from {\it Method of Fluxions}}
\label{fig:newton-colson-p60}
\end{figure}

\begin{figure} 
\vspace{6pt}
\centerline{
	\includegraphics[width=12cm]{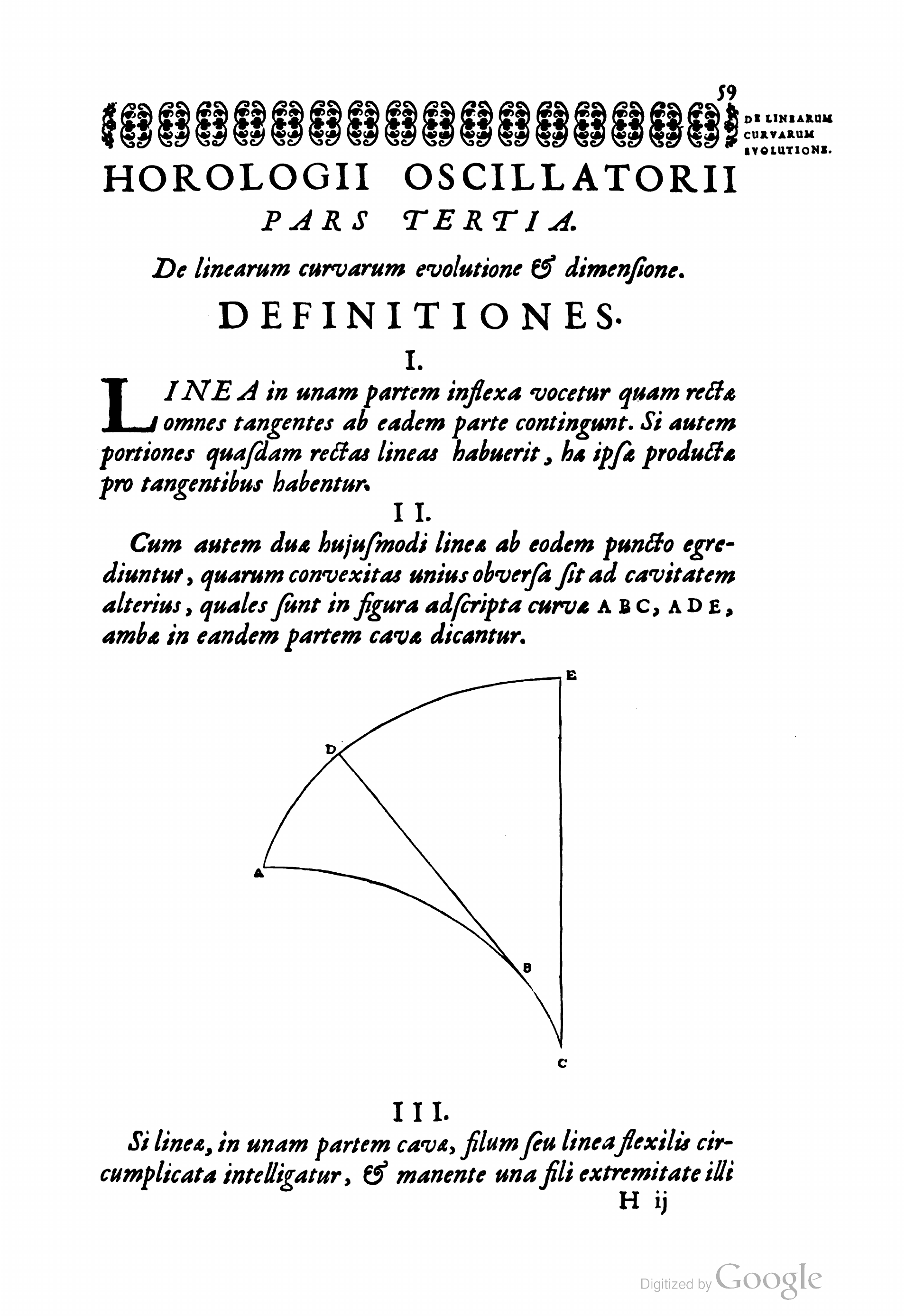}}
 \caption{Huygens center of curvature from {\it Horologium Oscillatorium}}
 \label{fig:huygens-p84}
\end{figure}

The first published account of curvature of a general curve was due to Huygens in 1673 (\cite{huygens1673}). In both Newton and Huygens the fundamental definition of the center of curvature (center of the osculating circle at a given point) is the intersection of normal lines to the curve near the given point on the curve (see the figures in Huygens p. 84 \cite{huygens1673} and Newton on p. 60 \cite{newton1736}, reproduced here in Figures \ref{fig:huygens-p84} and \ref{fig:newton-colson-p60}).

Huygens didn't have calculus at his disposal, but he made estimates in terms of normals at an approximating point (like the estimates of slopes of an approximating secant to a tangent line in differential calculus), and using these estimates he was able to compute the curvature for a variety of examples (cycloid, conic sections, etc.).  

\begin{figure} 
\vspace{6pt}
\centerline{
	\includegraphics[width=12cm]{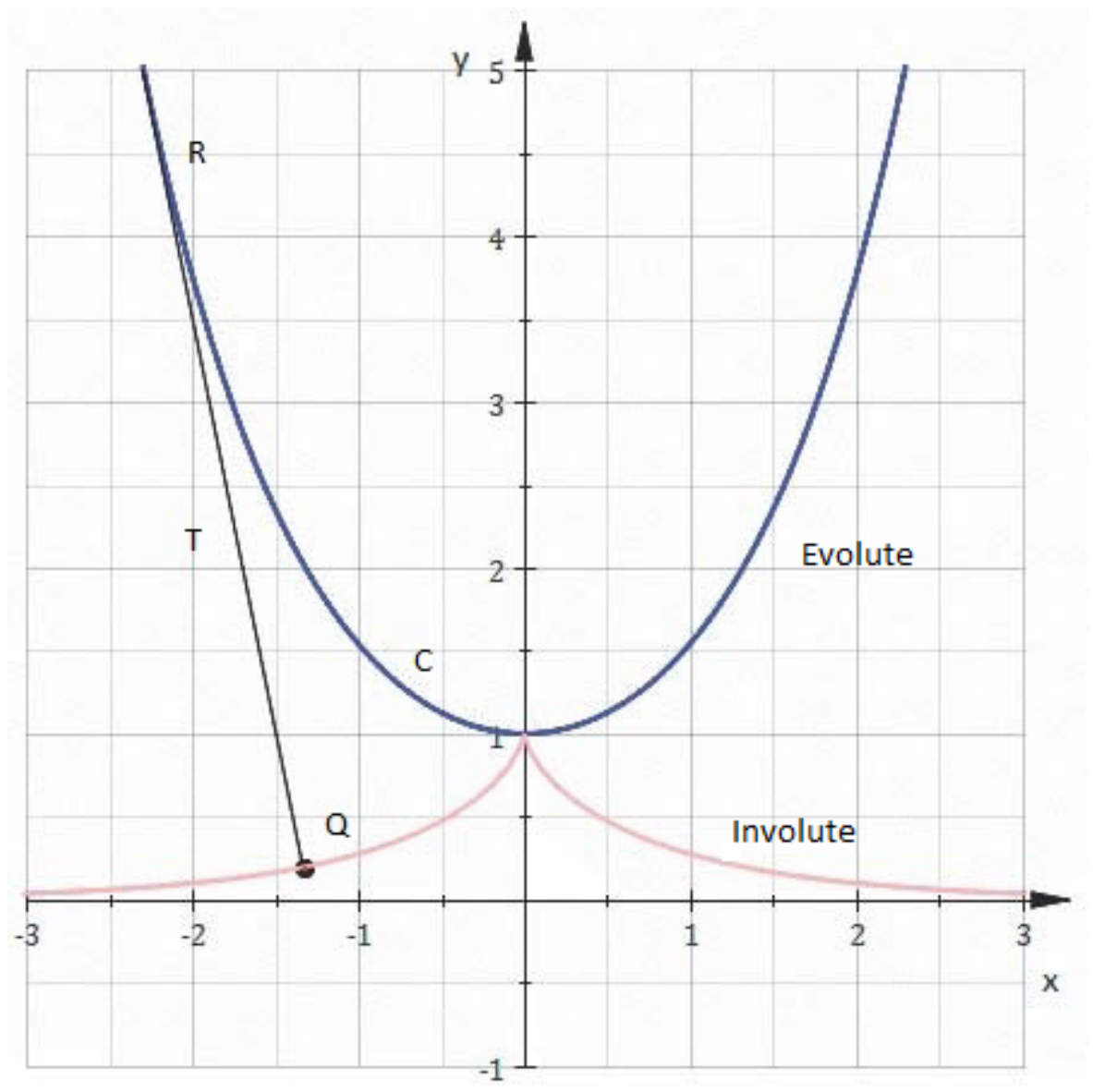}}
 \caption{An involute being generated by a string attached to the curve \(C\) (called the evolute)}
 \label{fig:involute}
\end{figure}

An interesting historical point is how Huygens came to study this phenomenon.  He had built some 16 years before the appearance of his monograph \cite{huygens1673} one of the most important clocks in history: a pendulum  whose motion is isochronous. That is the swing of the pendulum has a constant period of repetition. Huygens showed that a simple pendulum, whose pendant moves in a circular arc has a period that depends on the size of the oscillations, whereas if the pendant moves in the arc of a cycloid, then the period is fixed independent of the size of the oscillation.  The method Huygens used (which he patented in 1657) for making the pendant move in a cycloidal path was to have the path be the {\it involute} of a curved plate (which was also a cycloid), i.e., the curve traced out by a fixed string moving from a center attached to a given curve, where initially, the fixed string lies along the given curve, and moves away from it, with the free straight line portion of the string being continuously tangent to the given curve (see the illustration in Figure \ref{fig:involute}). 
  The curve \(C\) in Figure \ref{fig:involute} is called the {\it evolute} (which generates the involute  traced out by the point \(Q\)
by the motion of the string).  The problem Huygens posed and solved was; given the involute, find the evolute, i.e., find the generating curve.  Now the straight line \(T\) is normal to the involute at the point \(Q\) (as Huygens showed), and, at the point of contact at point \(R\), \(T\) is tangent to \(C\).  Thus \(T\) is normal to \(C\) at \(Q\) and \(R\) can be seen to be the intersections of the normals close to \(Q\) (as both Huygens and Newton showed).  Hence \(R\) is the center of curvature of the curve \(C\) at the point \(Q\), and the evolute \(C\) is the locus of centers of curvature of the involute at points near \(Q\).

\begin{figure} 
\vspace{6pt}
\centerline{
	\includegraphics[width=12cm]{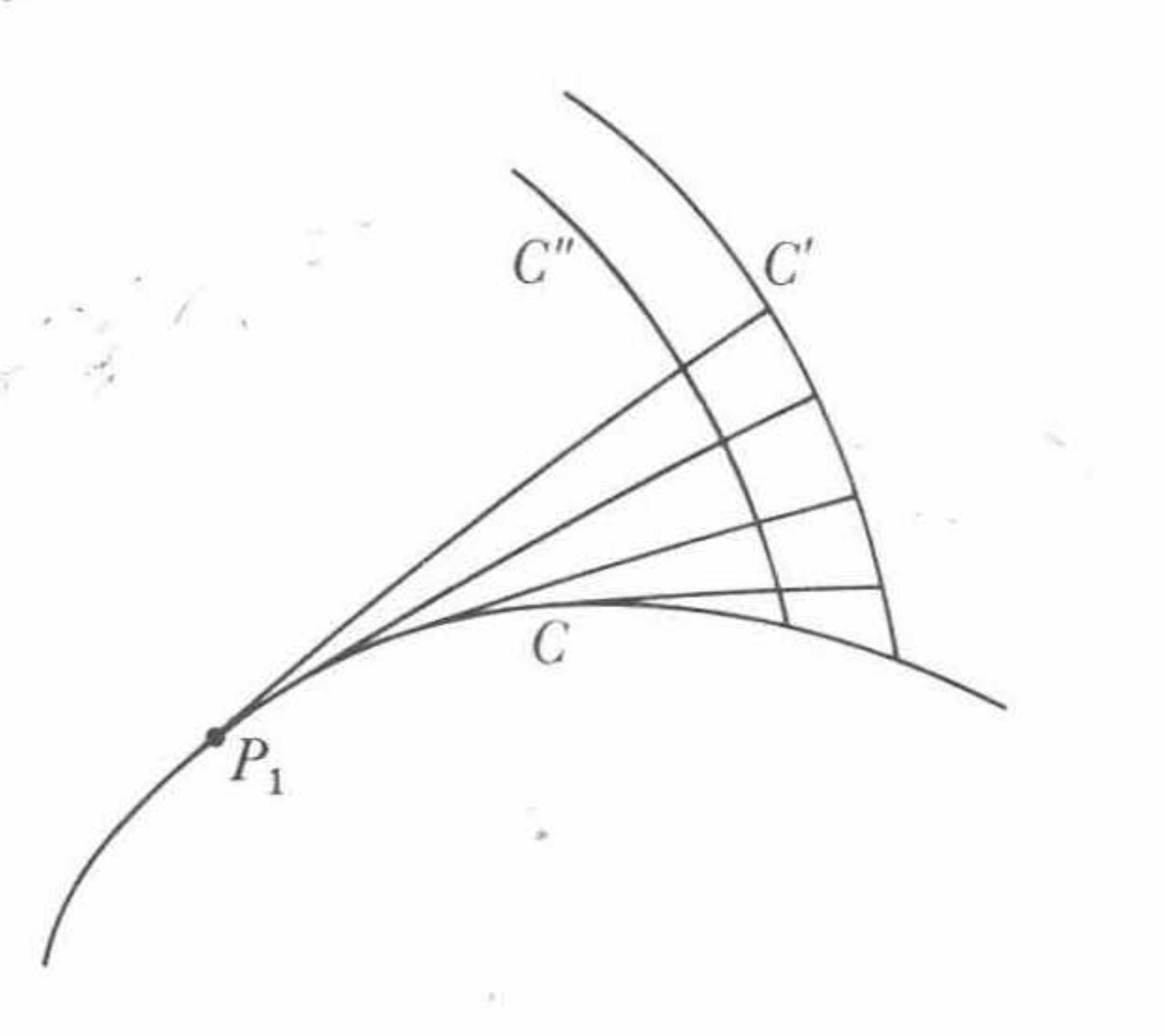}}
 \caption{Involutes are orthogonal to the generating string)}
 \label{fig:klinep555fig}
\end{figure}

\begin{figure} 
\vspace{6pt}
\centerline{
	\includegraphics[width=12cm]{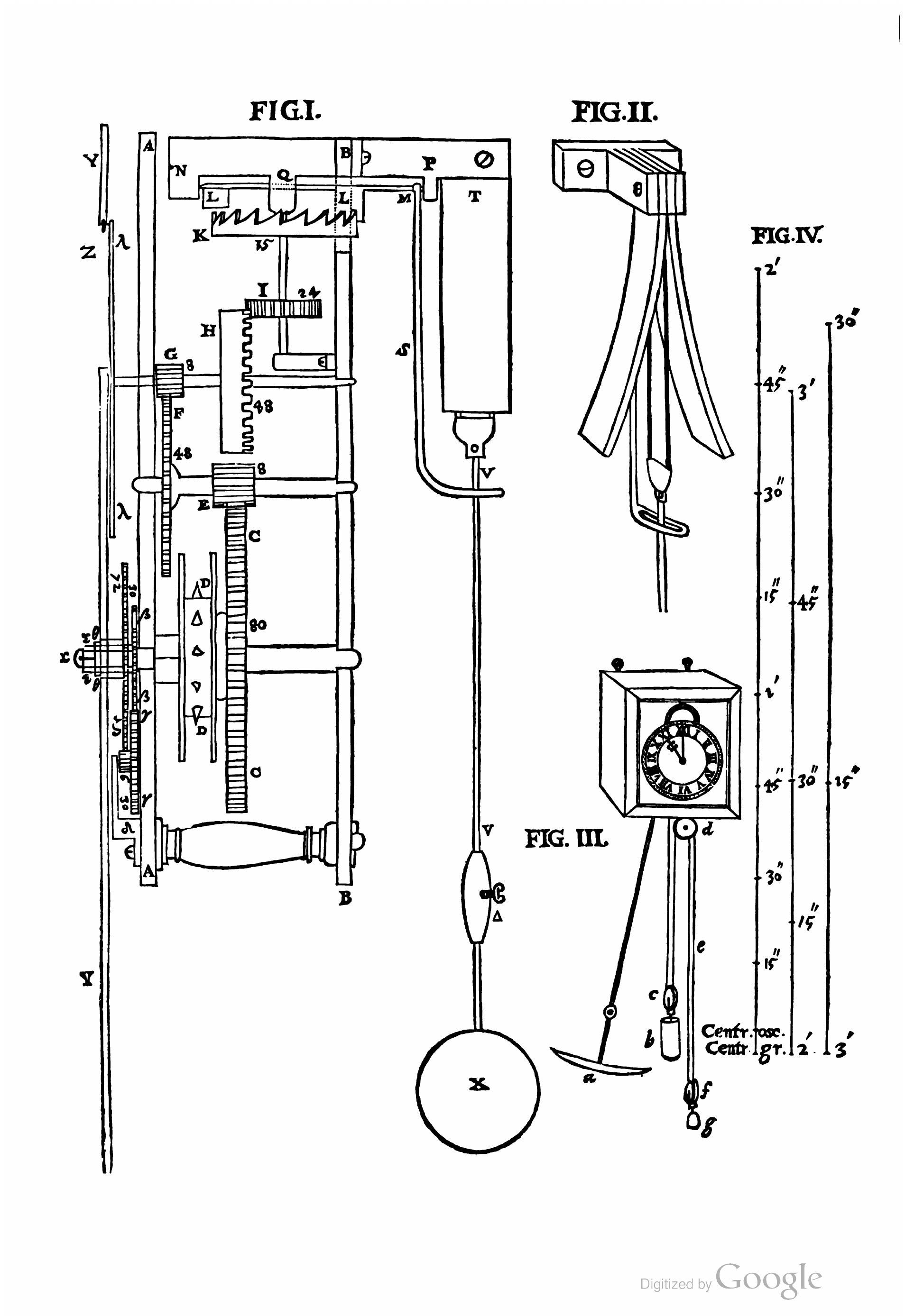}}
 \caption{Page 4 of Huygens' book {\it Horologium Oscillatorium}}
 \label{fig:huygensp4}
\end{figure}

In the second illustration of an involute in Figure \ref{fig:klinep555fig} one sees two ``parallel" evolutes, the curves \(C'\) and \(C''\) being generated from the curve \(C\), and one can see that the evolutes are orthogonal to the generating string at the intersection points (as was proved by Huygens). 
Looking at the illustration from p4 (Figure \ref{fig:huygensp4}) of Huygens' book \cite{huygens1673}
one sees in Figure II of the diagrams in Figure \ref{fig:huygensp4}  the cycloid-shaped curve from which the pendant of the pendulum sweeps out the involute which is the cycloidal motion of the pendant. 

Huygens calculated the evolutes for a number of examples, independent of his specific example he used in the design of his clock.  

Some 2000 years earlier, Apollonius
 was able to compute the curvature of the classical conic sections in Chapter V of his famous monograph on conic sections \cite{apollonius1896}.   Apollonius was in fact trying to solve a different set of problems and curvature was not explicitly discussed.  Namely, in the Heath translation of Conic Sections \cite{apollonius1896}, Heath shows what Apollonius
 did in modern notation. More particularly on p. 171 one finds that for the parabola of the form
\[
\frac{1}{2a}y^2=x
\]
the evolute (locus of centers of curvature) of this parabola has the form:
\[
27ay^2=4(x-2a)^3,
\]
 which is a semicubical parabola.  He finds similar formulas for the ellipse and hyperbola.
Here Apollonius was  studying the behavior of normals to conic sections.  He showed that each conic section has a unique normal passing through each point.  He {\it defined} a normal as being a straight line which was either a local maximum or a local minimum length straight line from some point not on the curve.  He then showed that such a line was indeed perpendicular to the tangent line at the given point.  This leads, by an interesting argument, to the conclusion that Apollonius has calculated the points of the evolute, as Heath points out very explicitly.

\subsection{Curvature of curves in space}
\label{sec:curvature-space}
\begin{figure} 
\vspace{6pt}
\centerline{
	\includegraphics[width=12cm]{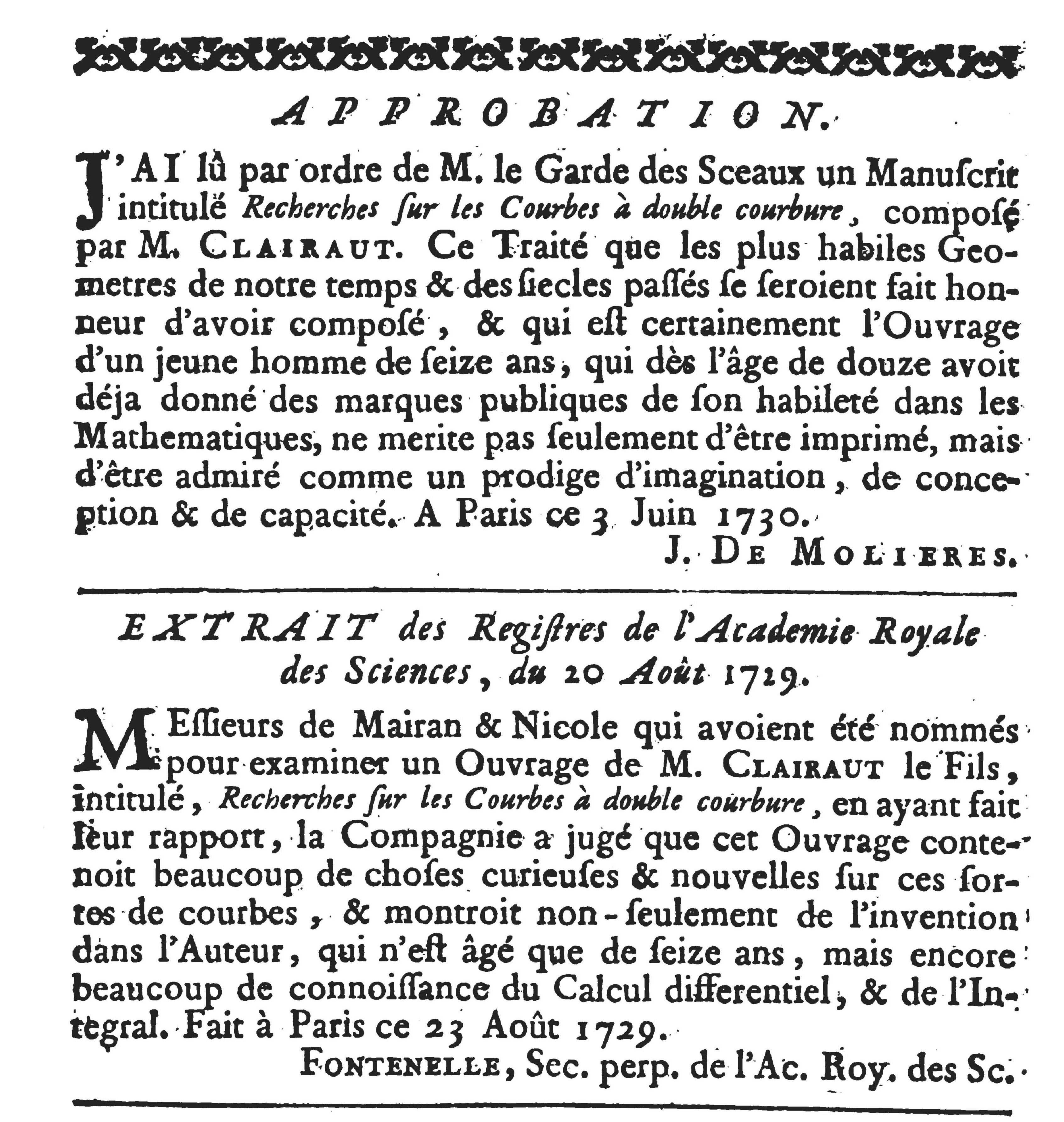}}
 \caption{Excerpt from beginning of Clairaut's book {\it Recherche sur les courbe \`a double courbure}}
 \label{fig:clairautp13par}
\end{figure}

From the time of Newton curvature of a curve in the plane became a standard object of mathematical investigation. The first step in investigating the differential geometry of curves in \(\BR^3\) was taken by Alexis Claude Clairaut (1713--1765) in his book {\it Recherche sur les courbe \`a double courbure} \cite{clairaut1731} written when he was only 16 years old and published two years later, following up on work he had started when he was 12 years old.  We know this from the ``Approbation" at the beginning of the book, written by two of the reviewers of the book, and the page where this appears, following the Preface, is the only place Clairaut's name appears in the book, not on the title page!  See Figure \ref{fig:clairautp13par}.
Clairaut called curves in \(\BR^3\) {\it curves with double curvature}, and he says in his book that he was inspired by Descartes, who suggested space curves could be studied in terms of their projections on two orthogonal planes.  Clairaut studied the tangent line to a curve, its arc length and the infinite variety of normal lines in the plane perpendicular to the tangent line.

The next steps in the study of space curves were taken by Euler, who primarily looked at space curves which were defined as the intersections of surfaces in \(\BR^3\) (see Volume 2 of Euler's {\it Introductio} of 1748 \cite{euler1748}).  Michel Ange Lancret (1774--1807) singled out in 1806 the three principal directions of a space curve at any point (tangent, normal, and binormal), and formulated the additional notion of torsion of a curve \cite{lancret1806}. The final steps in the study of space curves were taken by Augustin-Louis Cauchy (1789-1857) in 1826 in his {\it Le{\c{c}}ons sur les Applications du Calcul Infinit{\'{e}}simal a la G{\'{e}}om{\'{e}}trie} \cite{cauchy1826}, and by Serret \cite{serret1851} and Frenet \cite{frenet1852} in their back-to-back papers in 1851 and 1852. Cauchy gave us the formulation of space curves we use today (without the vector notation), and Serret and Frenet gave the final form to the structure equations (which today bear their name, the Frenet-Serret equations), which brought together the formal characterization of space curves in terms of the three principal directions of a curve and its curvature and torsion.

\subsection{Curvature of a surface in space: Euler in 1767}
\label{sec:curvature-surface}
\begin{figure} 
\vspace{6pt}
\centerline{
	\includegraphics[width=12cm]{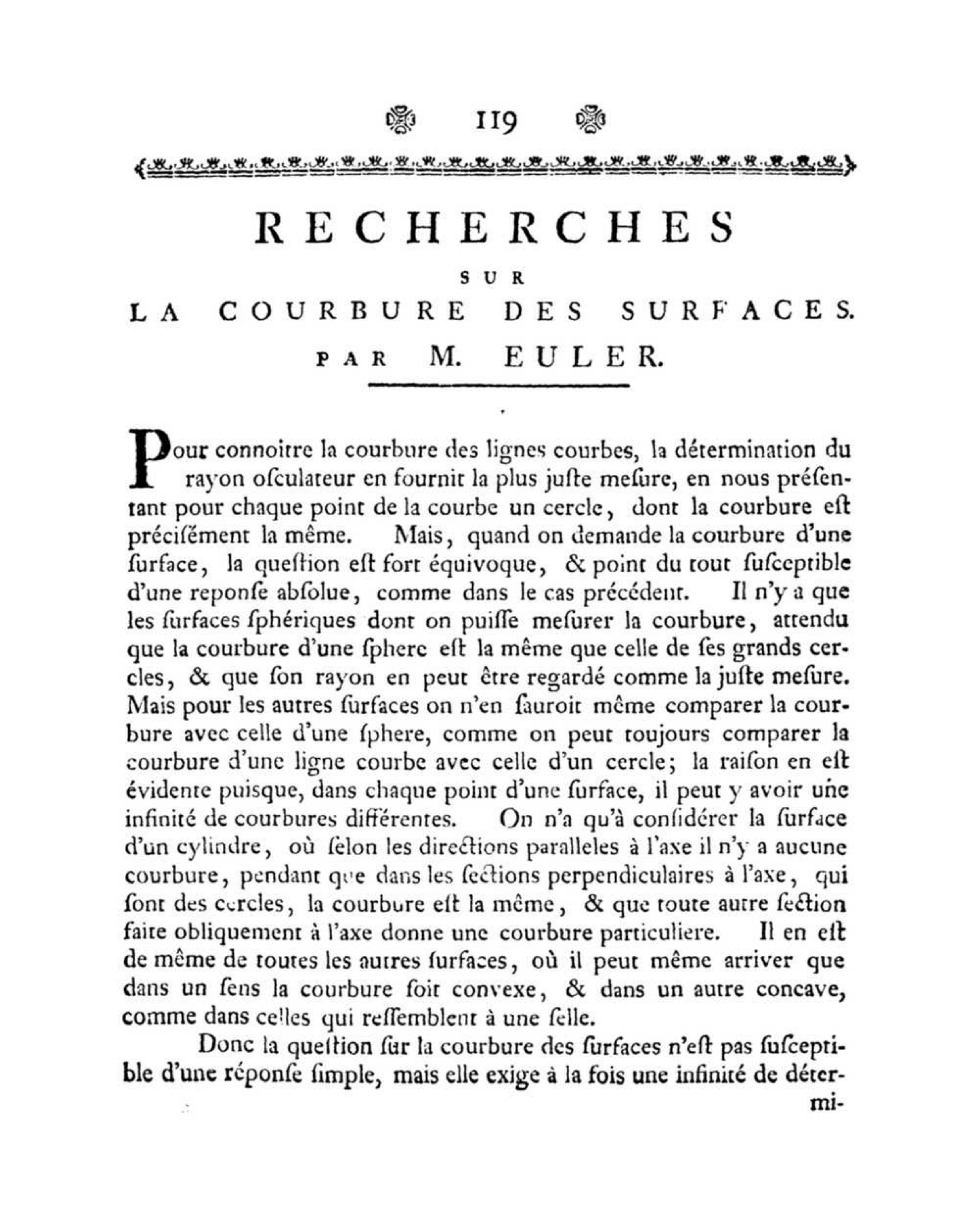}}
 \caption{The opening page of Euler's work on curvature}
 \label{fig:euler1767p119}
\end{figure}

This concept of curvature of a curve in \(\BR^3\) was well understood at the end of the eighteenth century, and the later work of Cauchy, Serret and Frenet completed this set of investigations begun by the young Clairaut a century earlier.  
The problem arose:  how can one define the curvature of a surface defined either locally or globally in \(\BR^3\)?   
An important contribution is made by Euler in his paper entitled ``Recherches sur la courbure des surfaces" \cite{euler1767} from 1767 (note this article is written in French, not like his earlier works, most of which were written in Latin).  Figure \ref{fig:euler1767p119} 
shows the first page of the article and we quote the translation here:
\begin{quote}
In order to know the curvature of a curve, the determination of the radius of the osculating circle furnishes us the best measure, where for each point of the curve we find a circle whose curvature is precisely the same. However, when one looks for the curvature of a surface, the question is very equivocal and not at all susceptible to an absolute response, as in the case above. There are only spherical surfaces where one would be able to measure the curvature, assuming the curvature of the sphere is the curvature of its great circles, and whose radius could be considered the appropriate measure. But for other surfaces one doesn't know even how to compare a surface with a sphere, as when one can always compare the curvature of a curve with that of a circle. The reason is evident, since at each point of a surface there are an infinite number of different curvatures. One has to only consider a cylinder, where along the directions parallel to the axis, there is no curvature, whereas in the directions perpendicular to the axis, which are circles, the curvatures are all the same, and all other obliques sections to the axis give a particular curvature. It's the same for all other surfaces, where it can happen that in one direction the curvature is convex, and in another it is concave, as in those resembling a saddle.
\end{quote}

In this paper Euler formulates quite clearly the problem of formulating a concept of curvature of a surface in \(\BR^3\). In particular, in the quote above one sees that Euler recognized the difficulties in defining curvature for a surface at any given point.  He does not resolve this issue in this paper, but he makes extensive calculations and several major contributions to the subject. He considers  a surface \(S\) in \(\BR^3\) defined as a graph
\[
z=f(x,y)
\]
near a given point \(P=(x_0,y_0,z_0)\).  At the point \(P\) he considers planes in \(\BR^3\) passing through the point \(P\) which intersect the surface in a curve in that given plane. For each such plane and corresponding curve he computes explicitly the curvature of the curve at the point \(P\) in terms of the given data.

He then restricts his attention to planes which are normal to the surface at \(P\) (planes containing the normal vector to the surface at \(P\)). There is a one-dimensional family of such planes \(E_\theta\), parametrized by an angle \(\theta\).  He computes explicitly the curvature of the intersections of \(E_\theta\) with \(S\) as a function of \(\theta\) , and observes that there is a maximum and minimum \(\k_\textrm{max}\) and \(\k_\textrm{min}\) of these curvatures at \(P\).%
\footnote{Moreover, he shows that the plane \(E_{\theta_\textrm{max}}\) is {\it perpendicular} to the plane \(E_{\theta_\textrm{min}}\) with the minimum curvature \(\theta_\textrm{min}\) (assuming the nondegenerate case; for a sphere all sections have the same curvature, and this statement would have no meaning).}
Moreover, he shows that if one knows the sectional curvature for angles \(\theta_1\ne\theta_2\ne\theta_3\), then the curvature at any given angle is a computable linear combination of these three curvatures in terms of the given geometric data.  This is as far as he goes.  He does {\it not} use this data to define the {\it  curvature} of the surface \(S\) at the point \(P\).  This step is taken by Gauss in a visionary and extremely important paper some 60 years later \cite{gauss1828}.

\section{Conclusion}
The ideas discussed above all play a major role in modern mathematics.  For instance, the classification of algebraic curves in the plane of degree two and three are important predecessors of what has become a major theme of contemporary geometry: to classify geometric objects. Moduli of Riemann surfaces, classification of topological and differentiable manifolds of various dimensions, including the famous Poincar\'{e} conjecture,  Thurston's classification conjecture for three-manifolds, the classification of two-dimensional complex manifolds (Kodaira--Spencer), and many other examples, are all instances of the classification of geometric objects. The work on differential geometry for curves and surfaces in three space that we describe here was an important prelude to the work of Gauss and Riemann on curvature of abstract manifolds in the 19th century, which has developed into the very rich field of differential geometry of the twentieth century.  We note that the recent solution of Grigori Perelman for the three-dimensional Poincar\'{e} conjecture, for instance, used the full power of differential geometry as a tool to solve this topological problem. The book by John Morgan and Gang Tian (Ricci Flow and the Poincar\'{e} Conjecture)  contains the details of Perelman's solution \cite{morgan-tian2007}.

\bibliography{references}
\bibliographystyle{plain}%

\end{document}